\documentclass[12pt,leqno,draft]{amsart} 
\theoremstyle{plain} 

\newtheorem{global-theorem}{Theorem}
\newtheorem{theorem}{Theorem}[section]
\newtheorem{lemma}[theorem]{Lemma}

\newtheorem{corollary}[theorem]{Corollary}

\newtheorem{prop-def}[theorem]{Proposition-Definition}
\newtheorem{lemma-def}[theorem]{Lemma-Definition}

\newcommand{\eop}{\ \hfill $\Box$}

\numberwithin{equation}{section}

\newcommand{\Xx}{{\mathcal X}}
\newcommand{\Cc}{{\mathcal C}}
\newcommand{\Ll}{{\mathcal L}}

\newcommand{\Ff}{{\mathcal F}}

\newcommand{\Uu}{{\mathcal U}}

\begin{document}

\author[C. Simpson]{Carlos Simpson}
\address{CNRS, Laboratoire J. A. Dieudonn\'e, UMR 6621
\\ Universit\'e de Nice-Sophia Antipolis\\
06108 Nice, Cedex 2, France}
\email{carlos@math.unice.fr}
\urladdr{http://math.unice.fr/$\sim$carlos/} 

\title[Gabriel-Zisman localization]{Explaining Gabriel-Zisman localization to the computer}


\maketitle

\tableofcontents

\section{Introduction} \label{sec-introduction}

We formalize the localization of categories, as in the book of Gabriel and Zisman \cite{gz},
with the Coq computer proof assistant. The purpose of this preprint is to provide some discussion of
this work. On the other hand, the computer files themselves
are attached to the companion preprint ``Files for Gabriel-Zisman localization''.
The text of that preprint consists mainly of the definitions and statements of results from the computer
files (in other words it is equal to the files with the proofs removed), plus some instructions for
compiling the files. 

There are several reasons for choosing this project. A certain amount of basic category theory was
done in the files attached to my previous paper on this subject \cite{fpccpag}. Thus it is natural to
look for some further topics to do in category theory. A long-range goal is to be able to do the theory
and practice of closed model categories. A glance at Quillen \cite{Quillen} suggests that the notion
of localization of categories {\em \`a la} Gabriel-Zisman is an important component of
the statements of some of Quillen's main results. Also in philosophical terms it is clear that Quillen
was influenced by Gabriel and Zisman, so it is reasonable to think that doing a computer formalization
of their construction of localization would be a good warm-up exercise.

A little bit of investigation into the bibliographical references for this construction has also turned up
another interesting reason to formalize it on the computer. It turns out that the full details of the 
construction (and specially of the calculus of fractions construction) have never really appeared in print.
Or at least, a search for these details  has not turned up any reference.
Of course it wouldn't be surprising to find a complete reference somewhere---you might say that this would be
the expected normalcy. Nonetheless it seems pretty clear that the vast majority of the very 
numerous mathematicians who use this theory every day 
haven't in fact read a text with the full details written down.

The notion of 
localization of a category is foundational for some of the most popular tools used by mathematicians
today: the homotopy category (of spaces, simplicial sets, or other things), and the derived category
(of an abelian category, coming in various flavors). It is surprising that the theory is so hard to find
written down in its integrality. This might contribute as part of an explanation for why the theories of
homotopy categories and derived categories are so much used and considered as ``black boxes''. 

One possible point of view would be to say that few have bothered to try to publish the full details
of the construction, because in a certain sense that just wouldn't be worth it: writing something down 
presupposes that there would be somebody interested in reading it; and writing down the full details of
an argument which is in essence straightforward, presupposes that a human reader would desire to, and be capable of, 
verifying in a meaningful way that the written text really did contain all of the details. 
Factors such as the total cost of publication also push towards leaving out much of this type of argument. 

In trying to write up the present note explaining the computer proof, it became evident that one had to agree
with the other published texts on this: the
full details of the argument just aren't sufficiently interesting to justify
the rather extensive linguistic effort which would be required to accurately convey them to a human 
reader, nor interesting enough for the reader to bother reading such an explanation. 
And this is with 
a fundamental piece of category theory more than $40$ years old. The arrival of the possibility
that the ``reader'' might be a computer changes this calculation. 
The computer is a perfect listener for an explanation that can be given as a sort of flow of
little arguments, sometimes with a necessary global strategy behind them, but always with lots of things
to remember, lots of referential notations to refer to various objects, and so forth. 

The purpose of the present preprint is to discuss our computer formulation, both of the general localization
construction and the special construction when there is a calculus of fractions. We don't pretend to give
all the details in the text---indeed we stop at about the same place as previous authors have. 
However, the details are necessarily all there in the computer files \cite{gzfiles}. 

\medskip

Historically, the notion of localization appeared informally in a somewhat different form in the
T\^ohoku paper of Grothendieck \cite{Grothendieck}
where he formalizes (to some extent) language which he attributes to Serre \cite{Serre} of 
working in a category ``modulo'' a subcategory.

After Gabriel-Zisman, the question of localization of categories has been 
treated in a number of references. Several people were
helpful in pointing out some of these in response to requests posted to the topology and category-theory
mailing lists. The references include books by L. and N. Popescu \cite{LNPopescu} \cite{PopescuAbCat}, 
H. Schubert \cite{Schubert}, and F. Borceux
\cite{Borceux}. Curiously enough the localization construction doesn't seem to appear
in \cite{MacLane}, although the underlying free and quotient category techniques are there. 
A classical reference which chronologically goes alongside Gabriel and Zisman is Verdier's thesis but which was only
recently published \cite{Verdier}. Verdier considers the case of localization of an additive category, inverting 
a multiplicative system which satisfies a two-sided calculus of fractions. The construction is similar
to the left-fraction construction.

The above list of references is undoubtedly partial.
It doesn't include very many of what are certainly
numerous research papers since the time
of \cite{gz} which may treat aspects of these issues in some detail (\cite{Benabou} is an example). 

Nonetheless, it is interesting to note the prevalence of formulations leaving ``to the reader''
parts of the proofs of details of the localization constructions. 
For example (the following all refer to the left or right fractions construction):
\newline
\cite{Verdier}, p. 117: ``La preuve est facile et laiss\'ee au lecteur qui pourra d\'emontrer de m\^eme 
la proposition ci-apr\`es \ldots ''. 
\newline
\cite{Schubert}, p. 260: ``Using (i), (ii), (iii) and (v), there is no difficulty in verifying
that the composite \ldots is well defined, \ldots ''.
\newline
\cite{PopescuAbCat}, p. 155: ``It is not difficult to see that with the equivalence relation (2) introduced
above one also has a well-defined composition law \ldots ''.

Another interesting reference is Pronk's paper on localization of $2$-categories
\cite{Pronk}, pointed out to me by I. Moerdijk. 
This paper constructs the localization of a $2$-category by a subset of $1$-morphisms
satisfying a generalization of the right fraction condition. In this case, there is no need to divide by an equivalence relation on
the set of $1$-arrows, because the appropriate arrows (i.e. pairs or what we call ``fraction symbols'')
are identified by the presence of $2$-cells making them equivalent. On the other hand,
knowing which $2$-cells are there would tell us which fraction-symbols need to be identified in the $1$-localization. Thus it seems likely that
from the high level of detail present in \cite{Pronk}, one could extract the complete set of necessary arguments for the localization
of a $1$-category. Nonetheless, the full set of details for the coherence relations on the level of $2$-cells is still too  much,
so the paper ends with: 
\newline
\cite{Pronk}, p. 302: ``It is left to the reader to verify that the above defined isomorphisms $a$, $l$ and $r$ are natural in
their arguments and satisfy the identity coherence axioms.''

We take the opportunity at this point in  the 
introduction to mention the colimit point of view about localization, even though it isn't
treated in our proof verification files.   
The contents of this discussion are touched upon by Gabriel-Zisman in 1.5.4 of Chapter I, see also 6.2 of 
Chapter II, and since then has become even more well-known. 

The localization $\Cc [\Sigma ^{-1}]$ can be viewed as a
pushout or colimit in the category of categories (one has to fix a universe $\Uu$ and consider
colimits in the category of $\Uu$-categories). To be precise, it is a pushout fitting into a
cocartesian diagram 
$$
\begin{array}{ccc}
\Sigma \times I & \rightarrow & \Cc \\
\downarrow & & \downarrow \\
\Sigma \times \overline{I} & \rightarrow & \Cc [\Sigma ^{-1}]
\end{array} 
$$
where here $\Sigma$ denotes the discrete set of maps to be inverted
(considered as a discrete category), $I$ denotes the category with objects 
$0$ and $1$ and one non-identity arrow $0\rightarrow 1$, and $I\subset \overline{I}$
denotes the completion to a category with two objects $0$ and $1$ joined by a single
isomorphism. In a certain sense the fact that this diagram is cocartesian just restates
the universal property of the localization. 

Ross Street in \cite{RossEmails} mentionned the above pushout point of view as well as
another closely related terminology, saying that the localization
is the ``coinverter'' of the 2-cell $\sigma$ which is the natural transformation from
$$
{\rm dom} : \Sigma \rightarrow C\;\;
\mbox{to}
\;\;
{\rm cod} : \Sigma \rightarrow C.
$$

The fact that the localization is a pushout implies that this operation is compatible with
colimits of categories. For example suppose 
$$
\begin{array}{ccc}
A& \rightarrow & B\\
\downarrow & & \downarrow \\
C& \rightarrow & P
\end{array} 
$$
is a cocartesian diagram of categories, and $\Sigma _A$, 
$\Sigma _B$, and $\Sigma _C$ are subsets of morphisms in $A$, $B$ and 
$C$ respectively such that $\Sigma _A$ maps into $\Sigma _B$ and $\Sigma _C$.
Let $\Sigma _P$ be the union of the images of $\Sigma _B$ and $\Sigma _C$.
Then the diagram 
$$
\begin{array}{ccc}
A[\Sigma _A^{-1}]& \rightarrow & B[\Sigma _B^{-1}]\\
\downarrow & & \downarrow \\
C[\Sigma _C^{-1}]& \rightarrow & P[\Sigma _P^{-1}]
\end{array} 
$$
is cocartesian. Let's stress that none of this is proven in the proof files.

An interesting case is when we take $\Sigma =Mor (C)$ to be the full set of morphisms of a category.
In this case the localization is the {\em groupoid-completion}, the universal groupoid with a functor
from $C$, denoted
$$
C^{\rm gr} := C[Mor(C)^{-1}].
$$
The compatibility with colimits stated above gives the following statement for groupoid completions:
if 
$$
\begin{array}{ccc}
A& \rightarrow & B\\
\downarrow & & \downarrow \\
C& \rightarrow & P
\end{array} 
$$
is a cocartesian diagram of categories then the diagram of groupoids
$$
\begin{array}{ccc}
A^{\rm gr}& \rightarrow & B^{\rm gr}\\
\downarrow & & \downarrow \\
C^{\rm gr}& \rightarrow & P^{\rm gr}
\end{array} 
$$
is cocartesian (either in the category of $\Uu$-categories, or in the category of $\Uu$-groupoids).
The groupoid completion $C^{\rm gr}$ is equivalent to the Poincar\'e fundamental groupoid of the
realization of the nerve of $C$ (denoted by $|C|$ for short).
In view of which, the above statement can be viewed as a sort of Van Kampen theorem
for fundamental groupoids in the style of R. Brown \cite{RBrown}. 
To make this view totally precise we would have to look at when the corresponding diagram 
of spaces $|A|$, $|B|$, $|C|$, $|P|$ is a pushout of spaces, which is sometimes but not always the case.

Ross Street also pointed out to me in \cite{RossEmails}, that the problem of localizing a noncommutative ring
is closely related (and somewhat similar to the Verdier case in that it brings in additive structure).
He sent me some notes treating in great detail the ring localization; it should be relatively easy to
transform those into a computer proof for this case. 
On the subject of notes, Clark Barwick mentions that he had written up some 
notes about Lemma 1.2 of Gabriel-Zisman; 
there are probably
(one would hope) a number of mathematicians who have done so.

Here is the plan of the paper. 
We will discuss first the general construction of the localization in a mathematical fashion; this
is followed by a section discussing the issues which arise in the computer formulation. Then
we come back to a mathematical discussion of the calculus-of-fractions construction, including discussion
of the subtleties which arise when we go towards the full details of the argument. The section after that
discusses the new issues which arise in the computer formulation, notably how we deal with the commutative
diagrams which one is tempted to use for the proof. In the last section, we mention very briefly
the contents of the remaining files in the present development.

\section{The general localization construction}

We recall in usual mathematical terms how the general construction of the localization of a category
works, taken directly from the first two pages after the introduction in Gabriel-Zisman \cite{gz}. 
Fix a category $\Cc$ and a subset $\Sigma \subset Mor (\Cc )$ of morphisms. We make no assumption
about $\Sigma$. We construct the localization, denoted $\Cc [\Sigma ^{-1}]$, as follows. 
Start by taking the disjoint union $Mor (\Cc ) \sqcup \Sigma$. The arrows in the first factor are thought
of as going in the forward direction, and the arrows in the second factor as going in the backward direction.
This allows us to create a directed graph whose vertices are the objects of $\Cc$ and whose set of edges is 
$Mor (\Cc ) \sqcup \Sigma$. Let $\Ff$ denote the free category over this graph. Recall that this means
that the objects of $\Ff$ are the vertices of the graph (thus, the same as the objects of $\Cc$),
and the morphisms of $\Ff$ are directed paths in the graph. 

Next introduce some relations on $\Ff$, and let $\Cc [\Sigma ^{-1}]$, which we denote as $\Ll$ for short,
be the quotient of $\Ff$ by these relations. The relations are trivial on the set of objects, that is to say
no different objects are put in relation and related morphisms share the same source and target. The relations
are introduced with the purpose of insuring the following properties for the quotient $\Ll$:
\newline
(1) \, the natural map on arrows $Mor (\Cc ) \rightarrow Mor(\Ff )$ (which is not itself a functor) projects
to a functor $\Cc \rightarrow \Ll$; and
\newline
(2)\, if $u\in \Sigma$ then the image in $\Cc$ of the backward edge corresponding to $u$ (that is,
the element of the second factor of the disjoint union) is inverse in $\Ll$ to the image of $u$ by the
functor in (1). 

The relations are chosen heuristically in a minimal way to accomplish this. In Gabriel-Zisman this process is
written in a compact way: the very process of taking the quotient category implies that we want the 
set-theoretical quotient of $Mor(\Ff )$ by the full set of relations, to be the set of morphisms of a category.
This in itself contains some properties of compatibility between the relations and the composition and 
identity operations. Call a relation which satisfies these properties, a {\em categorical relation}.
We can start by specifying an arbitrary list of relations and then take the closure under this condition,
that is the smallest categorical relation containing our list of relations. Once we have decided to do this,
we can list the germinal relations as follows: condition (1) requires that we specify two types of relations,
one for compatibility of the functor with the composition, and one for the compatibility of the functor
with the identity; and condition (2) requires again two types of relations, one each for the left and right
inverse properties. Thus we need to impose $4$ families of relations to start with. These are listed on
page 6 of Gabriel-Zisman (as well as in our proof file, see below). Given this list of relations, 
the quotient $\Ll$ is constructed by first completing to the smallest categorical relation $\sim$ containing the list,
then taking $Mor (\Ll ):= Mor(\Ff )/\sim $. It is straightforward to show that this defines a category
$\Ll$, with a functor $\Ff \rightarrow \Ll$ satisfying properties (1) and (2) above, since
we constructed the relation that way on purpose. We obtain a functor $P_\Sigma : \Cc \rightarrow \Ll$ sending
elements of $\Sigma$ to invertible morphisms in $\Ll$. 

So much for the construction. The next step is to state and prove an appropriate universal property.
On the first page of the construction they give the main properties which are (I am quoting):
\newline
``(i) $P_{\Sigma}$ makes the morphisms of $\Sigma$ invertible,
\newline
(ii) If a functor $F:\Cc \rightarrow \Xx$ makes the morphisms of $\Sigma$ invertible, there
exists one and only one functor $G:\Cc [\Sigma ^{-1}]\rightarrow \Xx$ such that $F=G\cdot P_{\Sigma}$.''

After explaining the construction in $14$ lines at the bottom of page 6 and the
top of page $7$ (which I have paraphrased above), 
Gabriel-Zisman jumped right up to a highly abstract formulation of the universal
property which we 
recopy here together with a few subsequent phrases:

\medskip

{\small 
``{\em 1.2.  Lemma:  For each category $X$, the functor $\underline{Hom}(P_{\Sigma},X): 
\underline{Hom}(\Cc [\Sigma ^{-1}], \Xx )\rightarrow 
\underline{Hom}(\Cc , \Xx )$ is an isomorphism from 
$\underline{Hom}(\Cc [\Sigma ^{-1}], \Xx )$ onto the full subcategory of $\underline{Hom}(\Cc , \Xx )$ 
whose objects are the functors $F:\Cc \rightarrow \Xx$ which make all the morphisms of $\Sigma$ 
invertible.}

``The proof is left to the reader. This lemma states more precisely conditions (i) and (ii). From 
now on \ldots ''
}

\medskip

Afterwards they pass immediately to the discussion of motivating examples like when the 
multiplicative system comes from a pair of adjoint functors.

This text (which totals less than a full page) is quite interesting from the point of view of the
problem of formalizing mathematics on the computer. In a very short space the authors have 
indicated, without error and indeed giving all of the necessary information, a relatively complex
mathematical construction, together with a very abstract statement of the universal property 
satisfied by this construction. Starting with the information given here, it is a
straightforward (and mathematically uninteresting) exercise to fill in all of the required details.

Creating a proof document to be read by a computer proof assistant raises a certain number of
mild difficulties. As a test, I have tried to attain the exact statement of Lemma 1.2 given above.
Before getting to a discussion of some details of this process in the next section, it is
interesting to note here that the time it took to do this was about a month, and the resulting total 
size of the $3$ proof files involved (\verb}freecat.v}, \verb}qcat.v} and 
\verb}gzdef.v}) is about 10,000 lines.

\newpage

\section{The computer formulation}

In order to formalize the general localization construction, we use the Coq proof assistant
\cite{coq}. We place ourselves in an axiomatic environment which  implements classical Zermelo-Fraenkel
set theory, and also relates it to the type theory of Coq. Concretely, the proof files attached to
\cite{gzfiles} include identical copies of the proof files of my earlier
set theory and category theory developments \cite{fpccpag}. The set-theoretical part of this 
starts with a file {\tt axioms.v} containing all of the axioms we assume (that is, no subsequent files
use the {\tt Axiom} or {\tt Parameter} commands). These axioms are intended to implement ZFC within Coq.
However, we furnish no formal proof of the fact that they do indeed do that, in other words that together with
the Coq ``calculus of inductive constructions'' type system, they furnish a mathematical system which
is consistent within the context of the usual ZFC axioms. It would be good to have such a proof, but that
seems to be complicated (due to the complicated structure of Coq) and possibly nontrivial due to 
certain aspects of Coq's type system such as cumulativity between sorts {\tt Prop} and {\tt Type}.
One would have to prove Conjecture 2 of Miquel and Werner \cite{MiquelWerner}. This is left to the reader!

The next question which is left open is to convince oneself that the definitions and statements of lemmas
contained both in the present files as well as in the category theory files, accurately represent
what the mathematician means when he speaks of categories, functors and so forth. This again 
may contain some nontrivial aspects
and is left to the reader. The accompanying preprint ``Files for Gabriel-Zisman Localization'' \cite{gzfiles} is
intended to help with this task: the textual part of this preprint consists of the Coq files, with
all proofs taken out. There one can look directly at the definitions and statements of lemmas,
which are the only parts which need to be understood in order to verify the meaning of what is
being said. (However, this is done only for the files concerning localization; in principle it might
be a good idea to have the same thing for the set-theory and category-theory developments but that
would be lengthy.)

While speaking about these foundational questions for the computer formulation, it is important to 
note that one could undoubtedly use any of a number of other environments for treating this question.
For example, it should be possible to proceed based on Saibi's category theory contribution 
\cite{Saibi} where sets are replaced by ``setoids'' (types plus equivalence relations). This is particularly
so in that one major element of the localization construction is the notion of quotient category. It is likely that
a setoid approach would simplify certain aspects, at the price of introducing other complications elsewhere.
We don't venture to predict how economical that would be on the whole. One should of course also
envision doing this type of formalization within other proof assistants (the list of which is getting very
long and we don't attempt to reproduce it here, see \cite{ctpm} \cite{umi} \cite{notices}).

For the reader who, at this point, still feels that a computer formalization can add something to the
question of verifying the mathematics underlying the localization construction, we now consider some 
details.

\subsection{Category theory}

We only
treat small categories, i.e. ones whose objects and morphisms form sets. In this context any  
distinction between small and big categories would be made by refering to a Grothendieck universe (which 
would itself be a set). However, in the current development 
we don't treat the question of when the localization of a category which is
big but has small $hom$ sets 
with respect to a given universe, is again big but with small $hom$ sets with respect to that universe. 
Thus, we are always working with sets and no foundational acrobatics come into play.

Most major notational questions about categories have already been dealt with in the category theory files.
A category is a $5$-uple consisting of the set of objects, the set of morphisms, the graph of the
partially defined composition operation, the graph of the identity operation, and a fifth set which 
is destined to contain any eventual extra structure one might want to include. The positions in the
$5$-uple are indicated by character strings (i.e. the $5$-uple is a function whose domain 
is a set of $5$ specific character strings corresponding to the $5$ places).  
This schema isn't the most economical: the necessary data (excepting the last structure variable)
is contained in the graph of the composition morphism. The goal is rather to achieve some rudimentary
standardization of the procedure for considering mathematical objects.

One feature of the category-theory encoding which is worthwhile to recall here is that
the set of morphisms is supposed to contain only objects 
\verb}u} which themselves are triples containing a source,
a target, and a third indicative element. 
This property is written \verb}Arrow.like u}. 
Here as before, these triples are realized as functions whose domain is a set of
$3$ character strings. This allows us to consider \verb}source u} and 
\verb}target u} for an arrow \verb}u}, independantly of the
category for which \verb}u} is a morphism. 
Here we have an economy of notation which 
has been extremely useful throughout the category theory development. In our discussion below we will 
encounter several places where a certain modification of the ``obvious'' approach is made necessary
by the \verb}Arrow.like} hypothesis. These modifications are easy to do once we are aware of the
phenomenon---which is why I am devoting some space below to these explanations. 

Similar notational considerations hold for functors and natural transformations. We refer to \cite{fpccpag}
for further discussion of these issues.

\subsection{The free category}

The first step in our current files is \verb}freecat.v} where we construct the free category on a graph. 
Since a morphism in the free category is a path in the graph, we need to implement the notion of path. 
This touches on what G. Gonthier explained was an important piece of their work formalizing
the 4-color theorem \cite{Gonthier}. However, the approaches are not the same since we are much less
concerned with efficient computation on these objects and more concerned with their theoretical manipulation. 
The notion of path also appears in T. Hales' recent formalization of the Jordan curve theorem \cite{HalesJordan}. 

To implement a notion of ``path'', 
we obviously  need
a theory of ``uples'', which are implemented as functions whose domain is an interval of natural
numbers of the form $[0,\ldots , n-1]$ where $n$ is the length of the uple. We define the function
\verb}Uple.create} to create an uple of length $l$ from a function \verb}f:nat -> E}, and a function
\verb}component} to get back the $i$th element of an uple. We need the \verb}length} function
as well as \verb}concatenate}. 

An important lemma is \verb}uple_extensionality} which says that
two uples of the same length with the same elements are the same; this allows us to prove
associativity of concatenation. J.S. Moore said for his ACL2 system that the first thing you would want
to prove was associativity of concatenation. In that type of system, uples or lists are 
inductive objects and associativity is a statement
proved by recurrence on the length. In Gonthier's paper \cite{Gonthier},
the notion of path is defined structurally so that associativity is automatic
by term reduction and need not be mentionned as a lemma.
In our case, the technical tool used to simplify the proof of
associativity, and most other manipulations of our uples which are functions of natural numbers, 
is the \verb}omega} tactic. 
The usefulness of this tactic, developped
by  Cr\'egut \cite{omegaref} based on an algorithm of Pugh \cite{pugh}, was pointed out to me by Marco Maggesi. 
It dispatches easily any arithmetic statement involving the standard operations and inequalities
on natural numbers. In our proof files this powerful tactic is abbreviated as \verb}om} which could
be interpreted alternatively as a reference to Buddhism or the Marseille soccer team.
Finishing the \verb}Uple} module is the operation \verb}utack} which corresponds to concatenation with
an uple of length one. This specific case enters often later so we treat it specifically. 

The notion of graph is relatively easy to encode.
A graph is a pair consisting of a set of vertices, and a set of edges. 
The edges of a graph are also supposed to be arrows. This situation is simple enough to provide a good
example of our general notational procedure  which we can recopy here:
{\small
\begin{verbatim}
Definition Vertices := R (v_(r_(t_ DOT ))).
Definition Edges := R (e_(d_(g_ DOT))).
Definition vertices a := V Vertices a.
Definition edges a := V Edges a.
Definition create v e :=
denote Vertices v
(denote Edges e stop).
Definition Graph.like a := a = create (vertices a) (edges a).

Definition Graph.axioms a := Graph.like a &
(forall u, inc u (edges a) -> Arrow.like u) &
(forall u, inc u (edges a) -> inc (source u) (vertices a)) &
(forall u, inc u (edges a) -> inc (target u) (vertices a)).
\end{verbatim}
}
We don't do any
theory of graphs beyond just the definition. 

Next we look at the paths which will make up the morphisms of the free category on a graph.
These are arrows whose third term are uples; and furthermore the uples will eventually 
(in the definitions \verb}arrow_chain} and \verb}mor_freecat}) be supposed to
be sequences of composable arrows in the graph, starting from the source of the arrow and ending at the
target. This situation requires a certain amount of specific treatment, for example we define a
version \verb}segment} of the previous \verb}component} function (and \verb}seg_length} instead of \verb}length}).
In general terms, this type of definition contracting two or some other small number of functions which
often occur together, necessitating the transposition of all of the lemmas concerning the pieces,
occurs all over the place and seems to be a general phenomenon. The composition operation for the free category
is defined by using concatenation of the underlying uples.  We also define the identity (whose uple has length $0$)
and prove all the various things needed to obtain the category axioms. We then would like to consider
functors from the free category into another category. For this, we need to define the operation of
composing together a composable sequence of arrows in a category (the definition \verb}mor_chain} is very
much like \verb}arrow_chain}). In this way we can state a universal property of
the free category on a graph (see the results concerning the construction
\verb}free_functor}). To close out this discussion we also consider (in the results concerning 
\verb}free_nt}) natural transformations
between functors whose sources are the free category. This is significantly easier because a
natural transformation is a function on objects, and the objects of the free category are just the
vertices of the graph.

\subsection{Quotient categories}

After the free category, the other main element of the construction is the notion of quotient category
(\verb}qcat.v}).
We are in a situation which is significantly easier than the general case: our relation
has no effect on the objects. In other words, we have a relation on the set of morphisms of a category, 
such that two morphisms which are related already have the same source and target. It is convenient to
distinguish two separate notions, denoted \verb}(cat_rel a r)} and \verb}(cat_equiv_rel a r)}. The
first means that \verb}r} is an arbitrary relation on the morphisms of a category \verb}a}, respecting
source and target. The second means that \verb}r} is an equivalence relation and compatible with 
the composition of \verb}a}. One important example of a \verb}cat_equiv_rel} is \verb}(coarse a)}
which puts in relation any two morphisms with the same source and target. The existence of this maximal relation
allows us by intersection to define the smallest \verb}cat_equiv_rel} containing a given \verb}cat_rel} \verb}r}.
We call this construction \verb}(cer a r)} (here \verb}cer} stands for the 
``\verb}c}ategorical \verb}e}quivalence \verb}r}elation'' on the category \verb}a} generated by \verb}r}).

In order to construct the quotient category of \verb}a} by \verb}r}, 
we need a manipulation called \verb}arrow_class}. The reason for this
is that in our notion of category, the morphisms are supposed to be \verb}Arrow.like}, i.e. triples having a source,
a target and an arrow. Thus we can't just say that the set of morphisms of the quotient category is the
usual quotient (i.e. set of equivalence classes) of the set of morphisms by the relation. 
Thus we define \verb}(arrow_class r u)} to be the arrow with the same source and target as \verb}u},
but whose third element is the equivalence class of \verb}u} for the relation \verb}r}. Now the
set of morphisms of the quotient category will be the image of this construction as \verb}u} runs
through the morphisms of \verb}a}. The definition \verb}(is_quotient_arrow a r u)} formalizes
the statement that \verb}u} is in the image. We also need a construction \verb}(arrow_rep v)} 
going in the other direction (see Lemmas \verb}related_arrow_rep_arrow_class} and
\verb}arrow_class_arrow_rep} saying that the two constructions are inverse in the appropriate sense).
We then define \verb}quot_id} and \verb}quot_comp}, the operations which will become the 
identity and composition for the quotient category. As usual, before trying to construct the category
it is a good idea to prove all of the necessary properties for these constructions. Then when
we construct \verb}(quotient_cat a r)} we prove destruct-create lemmas 
\verb}comp_quotient_cat} and \verb}id_quotient_cat} saying that the identity
and composition are \verb}quot_id} and \verb}quot_comp}. The destruct-create lemmas
\verb}ob_quotient_cat} and \verb}mor_quotient_cat} are proven after \verb}quotient_cat_axioms}
because the properties \verb}ob} and \verb}mor} include the category axioms for their first variables.

The module \verb}Quotient_Functor} does similar things for defining a functor \verb}qfunctor}
to the quotient category, and a functor \verb}qdotted} from the quotient category. The latter
terminology is intended to suggest that \verb}qdotted} is the dotted arrow which is filled in
in the universal property of the quotient category. Thus if \verb}f} is a functor and \verb}r}
a categorical equivalence relation on the category \verb}source f} we get a functor \verb}qdotted r f}
such that
\begin{verbatim}
source (qdotted r f) = quotient_cat (source f) r
target (qdotted r f) = target f
fcompose (qdotted r f) (qprojection (source f) r) = f.
\end{verbatim}
The unicity statement for the universal property says that if \verb}f} is a functor with
\begin{verbatim}
source f = quotient_cat a r
\end{verbatim}
then
\begin{verbatim}
f = qdotted r (fcompose f (qprojection a r)).
\end{verbatim}
There are no particular difficulties encountered in these arguments 
beyond the kind we have already discussed above.

Also contained in the file \verb}qcat.v} is a module \verb}Ob_Iso_Functor} dedicated to studying the
following situation. We have a functor \verb}f} and a category \verb}a}. We study the
pullback morphism induced by \verb}f}, denoted \verb}(pull_morphism a f)},
from \verb}(functor_cat (target f) a)} to
\verb}(functor_cat (source f) a)}. 
Recall that these constructions come from the file on functor categories \verb}functor_cat.v} in the category-theory
development. The purpose of this module is to contribute to the proof of Gabriel-Zisman's Lemma 1.2.
In particular,
we will want to apply this to the case where \verb}f} is the
functor from a category to its localization. Thus we assume that \verb}f} is an isomorphism
on objects. We develop a criterion for when \verb}(pull_morphism a f)} is fully faithful
and injective on objects (see the definition 
\verb}iso_to_full_subcategory}), or equivalently that it induces an isomorphism from
\verb}(functor_cat (target f) a)} to a full subcategory of
\verb}(functor_cat (source f) a)}. The equivalence between these notions is shown in Lemma
\verb}iso_to_full_subcategory_interp}. 

Intervening in the statement of the criterion is the construction \verb}add_inverses a s}.
This is the subset of morphisms of \verb}a} which are either already in \verb}s}, or else are
inverses in \verb}a} to morphisms in \verb}s}. Our criterion, 
stated
at the end of this module in Lemma 
\linebreak
\verb}iso_to_full_subcategory_pull_morphism_criterion},
says that if a functor \verb}f} is an isomorphism on objects, and if 
\verb}(add_inverses (target f) (mor_image f))} generates the category \verb}(target f)}, then 
for any category \verb}a} the pullback functor \verb}pull_morphism a f} is an isomorphism
onto a full subcategory. We will use this criterion, applied to the functor from a category
to its localization, to obtain half of the statement of Lemma 1.2. 

Finishing out the file \verb}qcat.v} is a module 
\verb}Associating_Quotient} which substantially recopies much of the definition
of quotient category. The only difference is that we don't start with a category
but only with a structure like a category but which doesn't necessarily satisfy the
associativity or left and right identity axioms. The idea is that the equivalence relation
will enforce these axioms. This construction is not needed for the general construction of localization,
but it will be needed later for the construction of the category of fractions. It didn't seem necessary
to go back and redo the whole quotient construction with this generality in mind: it is easier to 
recopy the relevant parts and change them. This might result in a file which is longer than
necessary, but one should keep in mind that the variable we are trying to economize is the
energy necessary to produce (or understand) the collection of files, not their total length.

\subsection{Construction of the localization}

Recall that the construction of the localization starts by looking at the graph whose edge set is the
disjoint union of the morphisms of \verb}a} with the elements of \verb}s}. Since these two sets are
anything but disjoint, we need some additional notation to implement the disjoint union. To this end, 
the first thing one notices at the start of the file 
\verb}gzdef.v} is the introduction of two sets \verb}Forward} and \verb}Backward}.
These are character strings (which are elements of \verb}E} hence sets, see \verb}notation.v} \cite{settheory}).
An edge of
the graph is either a ``forward edge'' corresponding to a morphism in \verb}a}, or else a ``backward edge''
corresponding to an element of the localizing system \verb}s}. 
The obvious thing is to try putting
\begin{verbatim}
forward_arrow u := pair Forward u
backward_arrow u := pair Backward u.
\end{verbatim}
However, {\em this doesn't work}. The reason is that the elements of the set of edges of the graph are supposed
to be \verb}Arrow.like}. To remedy this problem, we set
\begin{verbatim}
forward_arrow u := Arrow.create (source u) (target u) (pair Forward u)
backward_arrow u := Arrow.create (target u) (source u) (pair Backward u).
\end{verbatim} 
Notice that the source and target are
interchanged in the function \verb}backward_arrow}. 
The function \verb}original_arrow}
yields back the arrow we started with. Now \verb}loc_edges a s} is the set of such edges,
i.e. the union of the images of \verb}forward_arrow} and  \verb}backward_arrow} respectively on 
the morphisms of \verb}a} and on \verb}s}. The union is disjoint because \verb}Forward} and
\verb}Backward} are distinct.
Define the graph \verb}gz_graph a s} whose vertices are the objects of \verb}a} and whose edges are
\verb}loc_edges a s}.

From here, the construction basically follows the ordinary one, and doesn't really take up too much space. 
The free category on \verb}gz_graph a s} is called \verb}gz_freecat a s}.  The definition \verb}gz_rel a s}
is where the defining relations for the construction of the localization are listed. We recopy here
a lemma which rewrites that definition in a slightly more readable fashion. 
{\small
\begin{verbatim}
Lemma related_gz_rel : forall a s e f, localizing_system a s ->
related (gz_rel a s) e f =
((exists x, (ob a x & 
e = (forward_edge (id a x)) & 
f = (freecat_id x))) \/
(exists q, (inc q s & 
e = (freecat_comp (forward_edge q) (backward_edge q)) &
f = (freecat_id (target q)))) \/
(exists q, (inc q s & 
e = (freecat_comp (backward_edge q) (forward_edge q)) & 
f = (freecat_id (source q)))) \/
(exists u, exists v, (mor a u & mor a v & source u = target v &
e = (freecat_comp (forward_edge u) (forward_edge v)) & 
f = (forward_edge (comp a u v))))).
\end{verbatim}
}
The size of this text is
comparable to the size of the paragraph of \cite{gz} where the relations are listed.  Then \verb}gz_cer a s}
is the associated categorical equivalence relation, and \verb}gz_loc a s} is the quotient category. 
The functor from \verb}a} to \verb}gz_loc a s} is called \verb}gz_proj a s}.

The module \verb}GZ_Thm} is where we prove Gabriel-Zisman's Lemma 1.2. 
From the file \verb}qcat.v}, the
module \verb}Ob_Iso_Functor} furnishes the results necessary to prove the
part of Lemma 1.2 which says that pullback is an isomorphism onto a full subcategory.
As pointed out above, one of the delicate
points is that the statement of Lemma 1.2 involves the pullback morphism \verb}pull_morphism} between functor
categories. Functor categories were treated in the category theory development \cite{fpccpag},
and the place where we made use of that theory was in the module \verb}Ob_Iso_Functor} 
so we don't actually encounter it too much anymore here.  
This conclusion is stated in the present file as \verb}iso_to_subcategory_pull_gz_proj}, a corollary of the fact that
\verb}gz_loc} is generated by adding available inverses to the morphism image of the
functor \verb}gz_proj}.

The main part of the work done in the present module is to prove
the versal part of the universal property. Furthermore, rather than just giving a proof we would like
to have some useful notation. We start by introducing this notation for the free category: the
operation \verb}fr_dotted} corresponds to filling in the dotted line in a diagram expressing the versality
of the universal property. Similarly \verb}qdotted} did the same thing in the file \verb}qcat.v},
and putting them together we get a construction called \verb}gz_dotted} which expresses versality in the
following way.
Given \verb}a}, \verb}s} and a functor \verb}f} with \verb}source f = a}, we say
\verb}loc_compatible a s f} if \verb}f} sends elements of \verb}s} to invertible morphisms in \verb}target f}.
Then \verb}gz_dotted a s f} is a functor with
\begin{verbatim}
source (gz_dotted a s f) = gz_loc a s
target (gz_dotted a s f) = target f
fcompose (gz_dotted a s f) (gz_proj a s) = f
\end{verbatim}
The last property here, which is Lemma \verb}fcompose_gz_dotted_gz_proj}, corresponds to the 
versality property (ii) of \cite{gz}, page 6.
The uniqueness property (i) on page \cite{gz} was our Lemma \verb}gz_proj_epimorphic}. These are actually
the properties which are the most useful in practice. 

Our version of Lemma 1.2 is
given by two statements,\verb}iso_to_subcategory_pull_gz_proj} as noted above, and
for identification of the full subcategory image of \verb}pull_gz_proj},
the lemma \verb}ob_image_pull_gz_proj}. 

It might eventually be useful to have a more concrete description of natural transformations between functors
starting from the localization, but apart from the fact that it is implicitly contained in the statement of Lemma 1.2,
we don't treat this further here.

\section{Calculus of left (or right) fractions}

When I gave a talk in Nice about the computer formulation of the general localization 
construction, Charles Walter suggested that it would be interesting to compare the
formalization of the general construction of localization, with what would have to be done
to construct the localization in the presence of the habitual calculus of fractions conditions.
With this motivation I set out a while later to formalize the fractions construction from
Chapter 2 of Gabriel-Zisman. 

Contrarily to the general construction, it turned out (in my own opinion at least) that
filling in the details of the left-fractions construction involved some nontrivial
(if easy) mathematical thought, and drawing lots of diagrams. We don't draw diagrams
in the computer formulation (that might someday be possible but it is beyond the reach
of most computer proof assistants for the moment). As a replacement, we set up definitions of
situations involving several arrows of a category, which correspond to the diagrams we would want to draw. 
This will be discussed in the next section. 

In the present section we go into some detail about the mathematics of the problem, 
which stems from the fact that Gabriel-Zisman state their fraction construction under
a somewhat weak collection of hypotheses about the localizing system. 

The dual notions of left and right calculus of fractions are intended to be analogues of
the notion of multiplicative system for a commutative ring, which as was well-known leads
to a description of the localization as a set of ``fractions''. In the case of categories,
one would like to represent elements of the localizations as ``fractions'' or diagrams
$$
x \stackrel{v}{\rightarrow} y 
\stackrel{t}{\leftarrow} z,
$$
where by convention the arrows going backward are supposed to be in $\Sigma$. This diagram is viewed
as representing the morphism $t^{-1}v$ of $\Cc [\Sigma ^{-1}]$ so it is called a {\em left fraction symbol}.
We would like to have a nice set of conditions guaranteeing first of all that every morphism of the
localization can be written as a fraction; and second guaranteeing that the equivalence relation on
formal symbols $(t,v)$ whose quotient the set of morphisms $t^{-1}v$ is easy to understand. 
This collection of conditions is the {\em calculus of left fractions}. There will be a dual notion
of {\em calculus of right fractions} obtained by conjugating everything with the `opposite' construction.
Aside from the problem of implementing this conjugation in the computer formulation, we will focus on the
left-fraction case.

The conditions for a calculus of left fractions are given on \cite{gz} page 12, (2.2 a,b,c,d).
For convenience we reproduce them here:
\newline
(a)\, $\Sigma$ contains the identity morphisms of all objects of $\Cc$;
\newline
(b)\, $\Sigma$ is closed under composition;
\newline
(c) \, If 
$$
X' \stackrel{s}{\leftarrow} X \stackrel{u}{\rightarrow} Y
$$ 
is a diagram with $s\in \Sigma$ then there exists a commutative square
$$
\begin{array}{ccc}
X & \rightarrow & Y \\
\downarrow & & \downarrow \\
X' & \rightarrow & Y'
\end{array}
$$
such that the right downward map in the square is
in $\Sigma$; and
\newline
(d) \, If $f$ and $g$ are two maps from $X$ to $Y$ such that there exists $s\in\Sigma $ with $fs=gs$,
then there is a morphism $t:Y\rightarrow Y'$ in $\Sigma$ such that $tf=tg$.

Condition (c) says that every right fraction symbol can be completed to a left fraction symbol (in the 
commutative square,
$ X'\rightarrow Y' \leftarrow Y$ is a left fraction), that is dividing by an element of $\Sigma$ on the right
can be changed to division on the left. Condition (d) says that equalization on the right can be
changed to equalization on the left.

Most notable about this definition, specially in light of common practice in more recent times, is
what is left out. It is natural to require the following condition, which we call {\em three for two}:
\newline
(e)\, if $X\stackrel{g}{\rightarrow}Y \stackrel{f}{\rightarrow} Z$ is a composable pair of morphisms,
then if any two of $f$, $g$ and $fg$ are in $\Sigma$, the third one is too.

\smallskip

In general we can define the {\em saturation} of a set of morphisms to be the set $\Sigma ^{\rm sat}$
of all morphisms in $\Cc$ which become invertible in $\Cc [\Sigma ^{-1}]$. It is clear from the 
universal property that the
functor 
$$
\Cc [\Sigma ^{-1}] \rightarrow \Cc [(\Sigma ^{\rm sat})^{-1}]
$$
is an isomorphism, that is a set $\Sigma$ and its saturation share the same localization. 
It is also clear that for any set of morphisms, the saturation satisfies conditions (a), (b) and (e).
Thus from a certain perspective there would be no loss of generality in requiring that our set of morphisms
satisfy condition (e). For example, Quillen will later incorporate this condition as an important part of his notion of
``closed model category''. 

Nonetheless, Gabriel-Zisman don't make this requirement (and indeed they don't even speak of the
three-for-two condition (e) near here in the text). When you start to look closely at the details it
becomes clear that stating and proving the construction of the left-fraction localization in the absence of
the three-for-two condition is a bit of a challenge, one which they happily ask the reader to meet 
almost without saying anything about it, just subtlely giving the correct definition of the equivalence
relation so as to make it work.

Throughout the discussion of the left-fraction condition---as was the case for the general construction 
too---Gabriel-Zisman make reference to the construction of the sets of morphisms as being a direct limit construction.
We ignore this aspect here: it isn't treated in the formal proof development and we don't discuss it
in the informal presentation either. In fact it goes beyond the concrete character of the construction and it 
isn't clear whether it represents a useable piece of information (although that doesn't mean that it isn't conceptually
important). 

We now get to the description of the equivalence relations. We define a preliminary set 
of formal symbols $(t,f)$ consisting 
of two arrows having the same target, the first of which is in $\Sigma$. A left fraction symbol $(t,f)$ is
drawn as a diagram
$$
x \stackrel{f}{\rightarrow} y 
\stackrel{t}{\leftarrow} z.
$$
We would like to define the
set of morphisms of the left-fraction category to be the quotient of this preliminary set by an equivalence
relation (\cite{gz}, the top of page 13). 
Before stating the relation, notice that the ``source'' of the formal symbol $(t,f)$ is the source of $f$,
whereas the ``target'' of $(t,f)$ is defined to be the source of $t$. 
We call the common target of $t$ and $f$ the {\em vertex} of the symbol. 
The equivalence relation will preserve source and target. 
Two symbols $(s,f)$ and $(t,g)$ are said to be {\em equivalent} if there are maps $a$ and $b$
such that the source of $a$ is the vertex of $(s,f)$ and the source of $b$ is the vertex of $(t,g)$,
and $af = bg$, $as = bt$, and furthermore $as = bt$ is in $\Sigma$. Note that these conditions
automatically say that the targets of $a$ and $b$ are the same. 
See the second diagram on page 13 of Gabriel-Zisman. 

We can think of these conditions
as giving a symbol $(as,af)=(bt,bg)$ which is ``beyond'' both $(s,f)$ and $(t,g)$, and indeed this notion
is what we use in the proof development. We say that a symbol $(r,u)$ is {\em beyond} $(s,f)$ if
there exists a morphism $a$ whose source is the vertex of $(s,f)$ and target the vertex of $(r,u)$, 
such that $r=as$ and $u=af$. In this case we say that the morphism $a$ is an {\em intermediary}
from $(s,f)$ to $(r,u)$. 

A natural impulse would be to ask that the morphism $a$ (or the morphisms $a$ and $b$ in the definition of the
equivalence relation) be in $\Sigma$. This would be automatic from the conditions that $s$ and $as$ are in $\Sigma$,
if we had the three-for-two condition (e). However, if we try to do the construction in the absence of (e), 
we shouldn't ask that the intermediary morphism $a$ be in $\Sigma$ 
because then the construction
wouldn't work. A counterexample is discussed in the file \verb}lfcx.v}. 

If we have condition (e), then the proof that this defines an equivalence relation is relatively
straightforward. Without it, things are somewhat more tricky. 
The details necessary to overcome this problem must be considered as subsumed in the phrase
``It follows from (a), (b), (c),
(d) that this defines an equivalence relation \ldots '' in the middle of page 13 \cite{gz}. 
We will now explain how to see that. 

The main difficulty lies in proving that the equivalence relation is transitive. This may be rewritten
in terms of the notion of ``beyond'' as trying to show\footnote{Curiously 
enough, this type of reasoning closely
resembles the notions of reduction and normalization for $\lambda$-calculus; it might be interesting
to explore the analogy.}
that if 
two different symbols $(r,u)$ and $(r',u')$ are both beyond $(s,f)$, then there is a symbol 
$(q,v)$ which is beyond both $(r,u)$ and $(r',u')$.
In this case we have morphisms $a$ and $a'$ serving as intermediaries between $(s,f)$ and
$(r,u)$ or $(r',u')$ respectively. We would like to complete $a$ and $a'$ to a commutative square.
For this we would hope to use condition (c), which requires one of the morphisms to be in 
$\Sigma$. If we had condition (e) then this would be OK; in general an additional step is necessary.

Say that $(r,u)$ is {\em under} $(s,f)$ if
there is a morphism $a$ intermediary from $(s,f)$ to $(r,u)$ such that $a\in \Sigma$.
Note that ``under'' implies ``beyond'' but not necessarily vice-versa. 
The main observation is the following lemma, which gives a sort of weak replacement
for the $3$ for $2$ property, and its corollaries.

\begin{lemma}
\label{weak3for2}
Suppose $s\in \Sigma$ and $a$ is a morphism composable with $s$, such that $r:=as$ is in 
$\Sigma$. Then there exists a morphism $b$ such that $ba\in \Sigma$. 
\end{lemma}
{\em Proof:} Consider the diagram 
$$
\cdot \stackrel{r}{\leftarrow} \cdot \stackrel{s}{\rightarrow} \cdot .
$$
It is a right-fraction symbol because $r\in \Sigma$.
By condition (c) it can be transformed into a left-fraction symbol:
there exist morphisms $x$ and $t$ with $t\in \Sigma$ and
$xr=ts$. We would like to factorize $t$ into a product $ba$, however we may need to go
farther yet using condition (d). Our morphism $a$ goes from the target of $s$ to the
target of $r$, and we have 
$$
xas = xr = ts.
$$
In particular, we have two morphisms $xa$ and $t$ with the same source and target, 
equalized on the right by $s\in \Sigma$. By condition (d) there is a morphism
$c\in \Sigma$ such that $cxa = ct$. Recall that $t\in \Sigma$ so $ct\in\Sigma$,
and we can set $b=cx$ to obtain the lemma.
\eop

In the proof files, the argument of Lemma \ref{weak3for2} is integrated into the
proof of Lemma \verb}exists_lf_under} as in the following corollary. 

\begin{corollary}
\label{beyondunder}
Suppose $(r,u)$ is beyond $(s,f)$. Then there is another left fraction symbol $(t,v)$
such that $(t,v)$ is beyond $(r,u)$ and under $(s,f)$. 
\end{corollary}
{\em Proof:} (see \verb}exists_lf_under} in the proof files). 
Let $a$ be the intermediary morphism going from $(s,f)$ to $(r,u)$. Recall that $r$ and 
$s$ are in $\Sigma$, and $r = a s$. The lemma says there is another morphism $b$ such
that $ba\in \Sigma$. Let $t=br=(ba) s$ and $v =bu=(ba) f$.
\eop

\begin{corollary}
\label{beyondnormalization}
If $(r,u)$ and $(t,v)$ are both beyond $(s,f)$ then there is a symbol $(q,w)$ which is
beyond $(r,u)$ and $(t,v)$. 
\end{corollary}
{\em Proof:} By the previous corollary (and transitivity of 
``beyond'' which is easy) we may assume that $(r,u)$ is under $(s,f)$. Then 
(and this part is Lemma \verb}exists_lf_further} in the proof files) applying
condition (c) to the intermediate morphisms we obtain intermediate morphisms going
from $(r,u)$ and $(t,v)$ to a single $(q,w)$.
\eop

Transitivity of the relation follows easily from Corollary \ref{beyondnormalization}.
See \verb}lf_equiv_trans} in the proof files. 

A well-thought out 
direct argument for transitivity (which doesn't occupy too much space) is given in Borceux \cite{Borceux}
Proposition 5.2.4. The essential information is reduced to a single diagram (Diagram 5.4, page 185)
containing $9$ objects and $13$ arrows. 

For the definition of the composition, Borceux writes (p. 185):

{\small
``\ldots Moreover this definition is independent of the choices of $f$, $s$, $g$, $t$, $h$, $r$.
This is lengthy but straightforward: the arguments are analogous to those for proving the transitivity
of the equivalence relation defined on the arrows. We leave those details to the reader as well as
the checking of the category axioms \ldots .''}

This analysis is basically sound: once one has gotten over the hurdle discussed above, which first
shows up at the proof of transitivity, the remainder of the argument necessary for checking
well-definedness of the composition, associativity and identity axioms and so forth, presents no
further difficulties.  Nonetheless, it might be the case that the simplified presentation of the
proof of the transitivity of the equivalence relation, could have as a consequence that checking the
facts about the composition law becomes more involved (we needed to use techniques similar
to those for transitivity, in the proof of well-definedness of the composition for example). 
Similarly in \cite{Schubert} and \cite{PopescuAbCat}, the composition representative is constructed
but well-definedness and associativity of the composition are not verified in detail.

For completeness, we describe here some of the main points. 
Suppose we are given two left-fraction symbols which are composable:
$$
x \stackrel{f}{\rightarrow} y 
\stackrel{t}{\leftarrow} z
\stackrel{g}{\rightarrow} u 
\stackrel{r}{\leftarrow} v.
$$
Then the middle arrows give a right-fraction symbol which we can fill in to a square with a left-fraction
symbol going in the other direction:
\begin{eqnarray}
\label{compdiagram}
\begin{array}{ccccc}
 & & z & & \\
 & \swarrow & & \searrow & \\
y & & & & u \\
 & \searrow & & \swarrow & \\
  & & z' & & 
\end{array}
\end{eqnarray}
which in turn fits into the previous collection to yield a composite left-fraction symbol 
$(g'\circ f, t' \circ r)$:
$$
x \stackrel{f}{\rightarrow} y 
\stackrel{g'}{\rightarrow} z'
\stackrel{t'}{\leftarrow} u 
\stackrel{r}{\leftarrow} v.
$$
In order to define the composition, we make a choice of fill-in square and set the composition equal to 
the composite symbol $(g'\circ f, t' \circ r)$.  This composition rule is not associative, nor does it
satisfy the left and right identity relations.  On the other hand, modulo the equivalence relation
established above, the composition will become associative and unitary. 

The first and main step is to show that the composition is well-defined modulo the equivalence relation.
This has two parts: first that if we make two different choices of fill-in square then the resulting composite
symbols are equivalent; and secondly if we choose different representatives for the symbols which are being
composed, then the composites are equivalent. These proofs make use of the same kind of arguments as
we have described above, invoking things like Corollary \ref{beyondunder} when necessary. 
The reader can by now imagine why no authors have attempted
to write down the full text of these proofs in a forum destined for human readers. Those who are interested
may refer directly to the proof files. 

Once the well-definedness is established, the associativity is significantly easier at least on a conceptual
level. It suffices to look at the following diagram:
\begin{eqnarray}
\label{assocboard}
\begin{array}{ccccccccccccc}
\cdot & & & & \cdot & & & & \cdot & & & & \cdot \\
 & \searrow & & \swarrow & & \searrow & & \swarrow & & \searrow  & & \swarrow & \\
 & & \cdot & & 1 & & \cdot & & 2 & & \cdot & & \\
 & & & \searrow & & \swarrow & & \searrow & & \swarrow & & & \\
 & & & & \cdot & & 3 & & \cdot & & & & \\
 & & & & & \searrow & & \swarrow & & & & & \\
 & & & & & & \cdot & & & & & & 
 \end{array} 
\end{eqnarray}
where, along the top, are the three left-fraction symbols we want to compose. Choose the top two fill-in squares
denoted $1$ and $2$ 
first, which gives the middle row of arrows; then choose the bottom fill-in square $3$. 
Now the two different associated products may be obtained as follows (here we use the invariance under
choice of fill-in square): 
\newline
---one is obtained by using square $1$ to compose the first two symbols; then the composite
rectangle $2 + 3$ is a fill-in square for multiplying this first composite with the rightmost symbol;
\newline
---the other is obtained by using square $2$ to compose the second two symbols; and the composite
rectangle $1+3$ is a fill-in square for multiplying the leftmost symbol with this first composite. 

Both methods give as result the left-fraction symbol obtained using the composites along the
bottom edges of the big diagram. Thus, with the choices made as described above, the composition becomes
associative ``on the nose''; and because of the invariance of choices up to equivalence, we get that
composition is associative up to equivalence when we make arbitrary choices for the fill-in squares. 

The left and right unit conditions are proved similarly. 

We obtain a category of ``left fractions''. Defining the functor from our original category into the
left fraction category, and proving that the images of elements of the localizing system are invertible,
involve again some lemmas of a similar nature, whose proofs basically consist of setting up the appropriate
diagrams and using good choices for the fill-in squares to define the compositions in question. 
For all of these things, we are in agreement with all of the authors found so far, that it isn't worthwhile
to write a mathematical text for these proofs. The proofs may be found directly in the proof files
attached to \cite{gzfiles}. 

The construction of the localization by left fractions can be considered as the statement of a nontrivial
theorem about any localization (for example, about the general localization constructed previously).

\begin{theorem}
\label{lfproperty}
Suppose $\Cc$ is a category and $\Sigma$ is a multiplicative system satisfying the left-fraction conditions
(a)--(d) above. Let $\Cc [\Sigma ^{-1}]$ be a localization. Then the morphisms of $\Cc [\Sigma ^{-1}]$ 
have the following description. Every morphism can be written as a composition $t^{-1}\circ u$ where
$u$ comes from $\Cc$ and $t$ comes from $\Sigma$ (and their targets coincide). If, for two such pairs
$t^{-1}\circ u = r^{-1}\circ v$, then there exist morphisms $a$ and 
$b$ in $\Cc$ such that: $a$ is composable with $u$ and $t$;  $b$ is composable with $v$ and $r$;
$au = bv$; and $at = br$ and this is in $\Sigma$.
\end{theorem}

The proof is that this description holds by definition for the left-fraction localization we are discussing
in the present section. Then the universal properties show that any two localizations are isomorphic,
so the same description holds in any other localization. This theorem
is treated in the file \verb}gzloc.v}
(it is only there that we treat the fact that different localizations are isomorphic). 
It might be interesting to try to prove this description directly for the general construction of the localization.
The left-fractions conditions imply fairly directly that morphisms in the localization can be written as 
simple products. However, to verify the statement about the equivalence relation seems difficult. 

As a conclusion to this section, it is interesting to note that the mathematics behind the fraction construction
is not one hundred percent straightforward, as was the case for the mathematics behind the general construction.
On the other hand, it is commonly believed that the ``calculus of fractions'' construction is much more concrete and
easy to understand. A possible reason for this is that mathematicians are very attached to considering the ``size'' 
of the mathematical objects which they manipulate,
rather than the size of the associated mathematical theories. Since the arrows in the fraction construction are paths of length two
whereas the arrows of the general construction are paths of arbitrary length, people prefer to think about
the fraction construction (for example D. Pronk generalized the fraction construction to the case of $2$-categories \cite{Pronk}
but didn't mention generalizing the general construction).  This tendancy is similar to the constructionist or intuitionist
philosophy: even while admitting reasoning based on less constructive arguments, mathematicians of all philosophies 
gravitate towards smaller and more constructive objects when they are available.

\section{Formalizing the left-fraction construction}

The formalization is contained in the file \verb}left_fractions.v}, where  
\verb}Left_Fractions} is the first module treating all of the essential constructions and properties.
It starts with what is by now a fairly standard kind of definition, \verb}lf_symbol f t} is
an object (\verb}Arrow.like}, in fact) containing the pair $(f,t)$ and corresponding to the
left-fraction symbols used in the informal discussion above. 
The construction \verb}lf_choice a s r g} represents a choice of fill-in left-fraction symbol
creating a commutative square whose upper sides are the right-fraction symbol $(r,g)$. 

A left-fraction symbol has an additional object besides its \verb}source} and \verb}target},
which we call \verb}lf_vertex u}. This is the common target of the two morphisms involved.
The operation \verb}lf_extend a s p u} corresponds to composing both arrows of the left-fraction
symbol $u$, with a morphism $p$ whose source is \verb}source p = lf_vertex u}.

The \verb}lf_extend} enters into the definitions of the notions \verb}lf_beyond} and 
\verb}lf_under} as defined in the previous section. In turn we define \verb}lf_equiv}
as existence of a common symbol which is \verb}lf_beyond} the two in question.
Then comes the main part of the proof 
which is Lemma \verb}lf_equiv_trans}. This proof is done as described in the previous section
(the division into sublemmas is slightly different from what is done informally above;
we have referred above to the corresponding places in the proof files).

The main thing I would like to talk about in this section is the method we use for representing
situations which, in informal argument, would be represented by a commutative diagram drawn in the
text. Consider for example the definition 
{\small 
\begin{verbatim}
Definition fills_in a s u v w :=
has_left_fractions a s &
is_lf_symbol a s u & is_lf_symbol a s v & is_lf_symbol a s w
& source u = target v & source w = lf_vertex v &
target w = lf_vertex u &
comp a (lf_forward w) (lf_backward v) =
comp a (lf_backward w) (lf_forward u).
\end{verbatim}
}
This represents the diagram \ref{compdiagram} we have drawn in the previous section for defining 
composition. The variables \verb}u}, \verb}v} and \verb}w} are the three left-fraction symbols 
occuring in the diagram (the first two on the upper row and the last one completing the bottom).

In a similar way, the definition \verb}assoc_board a s u v w x y z} represents the 
diagram \ref{assocboard} we have drawn above for the associativity of composition.
Another important pair of diagrams are \verb}lf_lean_to a s e f g h i j}
and \verb}closes_lf_lean_to a s e f g h i j k l}. These correspond to diagrams, vaguely
shaped like ``lean-to's'', which we haven't drawn
above (due mostly to my lack of \TeX -nique), and which enter into the proof that the composition
is well-defined up to equivalence. 

The idea in all of these cases is to make a definition involving all of the objects occurring in the
diagram, which corresponds to commutativity of the diagram plus all of the other basic information 
it is supposed to satisfy (for example saying what the elements are, and that the sources and targets match up). 
In the cases \verb}fills_in} and \verb}assoc_board} we have chosen to represent the elements as being
the left-fraction symbols (i.e. they are pairs of arrows) whereas in the diagrams 
\verb}lf_lean_to} and \verb}closes_lf_lean_to} the elements are morphisms of \verb}a}. In both cases the
first variable \verb}a} is the category and the second variable \verb}s} is the multiplicative system
we are considering. 

Given this way of manipulating diagrams, we can then state the main steps in the proof.
For example, the main step of the well-definedness of composition is 
{\small 
\begin{verbatim}
Lemma lf_lean_to_closure : forall a s e f g h i j,
lf_lean_to a s e f g h i j ->
(exists k, exists l,(closes_lf_lean_to a s e f g h i j k l)).
\end{verbatim}
}
This lemma is then used in Lemma \verb}weak_rep_lf_equiv} (this fact doesn't show up in the text as
recopied in the preprint \cite{gzfiles} because it is inside the proof---which shows a limitation to the
idea of just copying definitions and lemma statements as a simplified presentation).
 
Often the diagram definition will occur as a hypothesis of a lemma. This is particularly true of the
hypothesis \verb}assoc_board a s u v w x y z} which occurs as a hypothesis in a few different intermediate lemmas
before the main statement of Lemma \verb}make_comp_assoc_board}. The definitions \verb}ffb_symbol}
and \verb}fbb_symbol} correspond to the two operations of juxtaposing two fill-in squares to
get a rectangle which composes to a new fill-in square giving the outer compositions in the
associativity statement (in the previous section this was where we looked at rectangles denoted 
$1+3$ or $2+3$). 

As a general matter, writing mathematics for the computer requires that we go to a notation which is completely
precise. Different strategies might be used for trying to keep a lid on the length of such notation.
It is interesting to note that in situations such as the present one, precise 
definitions such as \verb}ffb_symbol a s y z} replace vague statements such as ``juxtapose the squares denoted
$1$ and $3$ in the above diagram''.  One major problem in both cases is the problem of refering to pieces of
the diagram in question. The fact that we have included the various pieces as variables in our diagram
definition \verb}assoc_board a s u v w x y z} means that we can give small variable letters to each piece.
Thus the \verb}y} and \verb}z} in \verb}ffb_symbol a s y z} refer to places in \verb}assoc_board}.
In a mathematical text this referentiation
operation becomes cumbersome when we start to manipulate large numbers of objects: we are led to 
circumlocutions like ``the arrow on the upper left'' and generally get lost in the meanders of referencing conventions
of natural language (which are not totally precise nor sufficiently powerful). 
Another tempting solution would be to create a distinct mathematical object for the whole diagram.
This is pretty much what we have done with the notion of \verb}lf_symbol} for example. The drawback is that
we then need long names for the component pieces of the diagrams. In the case of \verb}lf_symbol} these were
the functions \verb}lf_forward} and \verb}lf_backward}. This approach was called for in the case of \verb}lf_symbol}
because of the frequency and diversity of manipulations we have to do with these objects. On the other hand,
with big diagrams which occur basically only once or a few times, it seems better to avoid allocating 
specific long names to their component pieces, and let the components be variables in a propositional representation
of the diagram as a functional property. The observations in this paragraph are not intended as strict edicts
but rather as ideas for one possible way to approach the language problems which are posed by 
computer formulation. 

Rather than going on in detail about the remainder of the construction 
(the module \verb}Left_Fraction_Category} is where we use the \verb}Associating_Quotient}
module to actually construct the category of left fractions; then we need to establish
its universal properties and so forth), we close this section with an observation about proof
technique. The proofs of the main lemmas referred to above are often rather long, because they
involve manipulations of large amounts of information. Since steps such as rewriting tend to 
produce residual goals, one arrives at a situation where there is an impossibly large number
of residual goals to treat at the end of a proof. Furthermore the treatment of these goals tends to
be highly repetitive. In this situation it is essential to maintain a certain level  of discipline
in the following sense: when one comes upon a rewriting situation which generates an additional goal,
one must go back to the start of the proof and add in that goal as an \verb}Assert} statement.
This can be done without changing the  numbering of hypotheses (which would be painful to correct at 
each occurence of this phenomenon) by considering the \verb}Assert} statements as ``sublemmas''
in the proof and naming them as such. Thus, rather than writing 
\begin{verbatim}
Assert (...)
\end{verbatim}
it is better to write 
\begin{verbatim}
Assert (lemA : ...).
\end{verbatim}
The proof itself becomes a location where there are many sublemmas. One can even recopy the sublemma
texts with their proofs, from one proof to another (when the proof contexts are going to be similar in
both cases). Once the required sublemmas are there, a combined rewriting tactic (such as
our abbreviations \verb}rw} or \verb}wr}) which does a rewrite and then tries \verb}assumption}
and \verb}trivial} on the subgoals, gives a proof where the number of auxiliary subgoals is
reduced significantly enough to be only ``very annoying'' rather than ``impossibly painful''.

\section{Further formalizations}

After the file on the left fraction construction, we include a few more files in the present
development. In \verb}gzloc.v}, we start by going back to some general considerations about 
localization. Lemmas whose names finish with \verb}_recall} are statements meant to recall definitions
from earlier files, so as to make reading of this part a little bit more self-contained. 
The first main corollary of the section is that two different localizations are isomorphic.
Recall that \verb}are_finverse a b} means that two functors are inverse (which is to say that their compositions
are {\em equal} to the identity---hence they establish an isomorphism between their sources and targets).
Thus the lemma \verb}are_finverse_dotted_choice} says that any two localization functors are isomorphic
(refer to \verb}source_dotted_choice} and \verb}target_dotted_choice} as well as \verb}fcompose_dotted_choice}
to complete this statement). 

The next step is to investigate the relationship between all of our various definitions, and the opposite
category and opposite functor structures. The main point here is to define the calculus of right fractions,
by conjugating left fractions with ``opposite''. This is of course completely straightforward, but many statements
require lengthy proofs, which is certainly evidence that our overall setup is not really optimal. 

In the part of the \verb}gzloc.v} file starting with the definitions
\verb}lf_vee}, \verb}lf_vee_image} and \linebreak
\verb}lf_vee_equivalent}, we treat Theorem \ref{lfproperty}. As pointed out above, the proof is
based on the fact that our left-fraction construction of a localization is isomorphic to any other 
localization due to the universal property. The statement of
Theorem \ref{lfproperty} is contained in Lemmas \verb}left_fraction_description_for_loc} 
(for the case of an arbitrary localization)
and \verb}left_fraction_description_gz_proj} (for the case of the original general construction).
Dualizing (by applying the opposite construction discussed in the previous paragraph)
we obtain the corresponding results for right fractions.  

The last file of the localization discussion is \verb}lfcx.v}. Here we construct a little counterexample
to one of the technical points encountered in the left-fractions construction. This part of the formalisation
takes us back to one of the basic points of our approach to types and set-theory, namely that we integrate
the inductive creation of types in {\sc Coq} into our axiomatization of set theory. This allows us to manipulate
small finite sets by creating them as inductive objects. In this way we can do a relatively large amount of
case analysis by defining recursive tactics using the \verb}Ltac} tactic language. In this way many of the proofs
of properties of our constructed objects are very short lists of tactics (which take some time for the
computer to digest). This could perhaps be thought of as a very very baby version of Gonthier-Werner's techniques
which they have applied to the 4-color theorem. The conclusion of the file is existence of a category
\verb}a} and localizing system \verb}s}, satisfying left fractions, but with two left-fraction symbols \verb}u}
and \verb}v} with \verb}lf_beyond a s u v} (thus \verb}u} and \verb}v} project to equivalent arrows in the
localization), however there is no symbol \verb}w} such that \verb}lf_under a s u w} and \verb}lf_under a s v w}.
This shows that we need to use the condition \verb}lf_beyond} rather than \verb}lf_under} (see the discussion of this
point above). 

The last file of the development is \verb}infinite.v}. This is a complete digression from the subject: here
we prove the fundamental properties of cardinal arithmetic for infinite cardinals, namely that the union or product of
two infinite cardinals has the maximum of the two cardinalities. By Russell's paradox on the other hand,
the powerset of an infinite cardinal is strictly bigger. Along the way we do a certain number of basic properties
of finite and infinite sets, cardinals, and ordinals. It is beyond the scope of the present preprint to go into
further detail about the strategy of proof; and (contrarily to the case of
localization of categories as we have seen) this is a subject in which  there is no lack of different treatments
in the literature---which we have not tried at all
to index in the references.

\end{document}